\documentclass[12pt]{article}
\usepackage{amsmath,amsthm,amsfonts,amscd,amssymb,eucal,latexsym}
\begin{document} 

\def\C{\mathbb C} 
\def\N{\mathbb N} 
\newtheorem{thm}{Theorem}[section] \newtheorem{corol}[thm]{Corollary} 
\newtheorem{prop}[thm]{Proposition} \newtheorem{lem}[thm]{Lemma}
\theoremstyle{remark}
\newtheorem{rem}[thm]{Remark} 
\theoremstyle{remark}
\newtheorem{exa}[thm]{Example} 
\newtheorem{ack}{Acknowledgement} \renewcommand{\theack}{}
\title{Representations of the Direct Product of Matrix Algebras} 
\author{Daniele Guido\footnote{Supported by GNAFA and 
MURST. Present address: Dip. di Matematica, Univ. della Basilicata, 
Contrada Macchia Romana, I-85100 Potenza, Italy, e-mail: 
guido$@$unibas.it} \ and Lars Tuset\footnote{Supported by CEE. 
Present address: Faculty of Engineering, Oslo University College, 
Cort Adelers Gate 30, 0254 Oslo, Norway, e-mail: larst$@$pc.iu.hioslo.no}
\\Dipartimento di Matematica,\\ Universit\`a di Roma ``Tor 
Vergata'',\\ I--00133 Roma, Italy} 
\date{} 
\maketitle 
\begin{abstract} Suppose $B$ is the unital algebra consisting of the 
algebraic product of full matrix algebras over an index set $X$.  A 
bijection is set up between the equivalence classes of irreducible   
representations of $B$ as operators on a Banach space and the 
$\sigma$-complete ultrafilters on $X$,  Theorem \ref{prop1}.   
Therefore, if $X$ has less than measurable cardinality (e.g. 
accessible), the equivalence classes of the irreducible   
representations of $B$ are labeled by points of $X$, and all 
representations of $B$ are described, Theorem \ref{Th2}.  
\end{abstract} 
\vfill\eject 
\section{Introduction} 
Suppose $B=\prod_{x\in X}B(K_x )$ is an algebraic product of full 
matrix algebras $B(K_x )$ over an index set $X$.  Clearly $B$ is a 
unital algebra under pointwise operations.  Consider now a 
representation $\pi$ of $B$ on a Banach space.  The purpose of this 
paper is to characterize such representations up to equivalence, see 
Theorem \ref{Th2}.  A simple characterization is, however, only 
possible under the condition that $X$ is a set of less than measurable 
cardinality.  While this condition imposes no real restriction on $X$ 
for practical applications of our result, since every set which can be 
constructed (by the ordinary operations in set theory, such as unions 
and powers) is of less than measurable cardinality, it is remarkable 
that such a non trivial condition arises in this context.  Theorem 
\ref{Th2} says that every representation $\pi$ of $B$ (where the index 
set $X$ is assumed to be of less than measurable cardinality) is 
equivalent to a direct sum (including if necessary infinite 
multiplicity) of representations obtained by restricting $\pi$ to some 
of the factors $B(K_x )$ of $B$.  As a corollary (see Corollary 
\ref{projned}) to Theorem \ref{Th2}, we conclude that when $\pi$ is 
irreducible, it is equivalent to the representation of $B$ obtained by 
projecting down on one of the factors of $B$.  We give a separate 
proof of this corollary as a shortcut to the reader.  We note that our 
results do not depend crucially on the fact that we represent 
unbounded elements with bounded operators, cf.  Remark 
\ref{rem:unbdd}.

The crucial notion in these investigations is that of a 
$\sigma$-complete ultrafilter.  To explain how these filters occur, 
let $\pi$ be an irreducible representation of $B$ on a Banach space 
$K$.  Define $F_{\pi}$ to be the collection of subsets of $X$ given by 
$F_{\pi}=\{ U\in X\, |\, \pi(\chi_U)=I_{\pi}\, \}.$
Here $I_{\pi}$ is the unit of $B(K)$ and $\chi_U$ is the element of 
$B$ defined by $\chi_U (x)=I_x$ for $x\in U$ and $\chi_U (x)=0$ for 
$x\notin U$, where $I_x$ is the unit of $B(K_x )$.  We show (Lemma 
\ref{ultra}) that $F_{\pi}$ is a $\sigma$-complete ultrafilter over 
$X$.  In fact, in Theorem \ref{prop1} we establish our second main 
result, which states that the assignment $\pi\mapsto F_{\pi}$ induces 
a bijection between the equivalence classes of irreducible 
representations of $B$ on Banach spaces and $\sigma$-complete 
ultrafilters over $X$.  This result is true for an arbitrary index set 
$X$.  However, it is known (see Lemma \ref{lem2a} and the paragraph 
prior to this lemma) that whenever $X$ is of less than measurable 
cardinality, then every such filter has to be principal.  This way an 
irreducible representation $\pi$ of $B$ singles out a point in $X$ 
(over which the maximal principal filter is based), and the 
representation obtained by projecting down on the factor of $B$ 
corresponding to this point, is equivalent to $\pi$.

We remark that in a very particular case our results were already 
known, namely when the spaces $K_x$ are one-dimensional.  In this case 
B is abelian, hence its representations are 1-dimensional.  Therefore 
our results imply that all multiplicative linear functionals on 
${\C}(X)$ are in $1:1$ correspondence with $\sigma$-complete 
ultrafilters on $X$.  As a consequence, all multiplicative linear 
functionals are given by Dirac measures if and only if $X$ has less 
than measurable cardinality.  Such result is contained in \cite{BBZ}, 
where it was generalised in a direction different from ours: the 
fibers of the direct product were allowed to be arbitrary algebras on 
an arbitrary field (with some assuptions on the cardinality of the 
field) but only one-dimensional representations, namely multiplicative 
linear functionals, were studied.  \footnote{Indeed we became aware of 
\cite{BBZ} when this paper was already finished.}

Even this particular case cannot be easily guessed by analogous 
results.  For instance, the set of characters of the algebra of 
complex valued functions on $X$ with finite support (or vanishing at 
infinity) can be identified with the points of $X$.  This is a 
particular case of the result that every character on the 
$C^*$-algebra of continuous complex valued functions vanishing at 
infinity on a locally compact Hausdorff space is a Dirac measure, cf.  
\cite{GM}.  If we instead consider the algebra of bounded complex 
valued functions on $X$, then the set of characters is the Stone-Cech 
compactification of the space $X$ equipped with discrete topology, cf.  
\cite{Wi}.  Thus in this case the characters are Dirac measures only 
when $X$ is finite.  It is remarkable that passing from bounded to 
unbounded functions amounts to passing from finite to less then 
measurable sets.

The idea of interpreting noncommutative algebras as noncommutative or 
quantum spaces has proved fruitful.  We only mention the theory of 
noncommutative differential geometry developed by A. Connes (cf.  
\cite{AC}), and the theory of quantum groups (cf.  \cite{Wo}).  In 
view of this one may conceive Corollary \ref{projned} as a first 
humble step towards \lq noncommutative set theory\rq\, with the \lq 
points\rq\, played by irreducible representations of the 
noncommutative product algebra $B$.  Finally, notice that such product 
algebras appear in the theory of discrete quantum groups, cf.  
\cite{tannakawor}.  The index set is then the equivalence classes of a 
complete family of the dual compact quantum group, or if you like, the 
equivalence classes of irreducible objects in a concrete tensor 
$C^*$-category with conjugates, see \cite{Y}.

\section{Filters and irreducible Representations} 
\label{sec2} 
Let $X$ be an arbitrary set, $x\to K_{x}$ a map associating with any  
point $x\in X$ a Banach space $K_{x}$, and let $B(K_x )$ be the 
Banach  algebra of bounded operators on $K_x$.   
Denote by $B$ the (unrestricted) algebraic direct product  
$B=\prod_{x\in X}B(K_x )$, i.e., $b\in B$ if $b$ is a function  
$b:X\rightarrow \cup_{x\in X}B(K_x )$ such that $b(x)\in B(K_x )$ 
for  all $x\in X$.  Clearly $B$ is a unital algebra with pointwise  
operations.  Whenever the Banach spaces $K_x$ are in addition 
Hilbert  spaces, the algebra $B$ is a $*$-algebra.  When all the 
Banach spaces  $K_x$ are the complex numbers ${\C}$, we denote the 
commutative  unital algebra $B$ by ${\C}(X)$, being all the 
complex valued  
functions on the set $X$. 
A representation $\pi$ of $B$ on a Banach space $K$ is a unital  
homomorphism $\pi :B\rightarrow B(K)$.  Two representations are  
equivalent if there is an invertible intertwiner between them.  A  
representation is said to be irreducible if the only idempotents  
intertwining it are zero or the identity. 
If $K$ is a Hilbert space and $\pi (B)$ is a $*$-algebra (e.g. when  
$B$ is a $*$-algebra and $\pi$ is a $*$-representation), then the  
commutator $\pi (B)'$ of $\pi (B)$ in $B(K)$ is a von Neumann 
algebra.   Hence it is generated by its (self-adjoint) projections.  
Therefore  irreducibility of $\pi$ is equivalent to the triviality of 
all (self-adjoint) projections intertwining it.  However, when $\pi 
(B)$ is not  a $*$-algebra (and $K$ is a Hilbert space), the two 
notions of  irreducibility do not coincide as the following example 
shows.   
\begin{exa} Consider $B={\C}^2$.  Denote by $<a>$ the  
universal unital complex algebra generated by the element $a$  
satisfying the relation $a^2 =I$, where $I$ is the identity.  Denote  
by $<b>$ the unital universal complex algebra generated by the 
element  $b$ satisfying the relation $b^2 =b$.  It is easy to see 
that  $b\mapsto a=I-2b$ extends to an isomorphism between $<b>$ and 
$<a>$.   Now every element of $<b>$ is of the form $\alpha I+\beta b$ 
for  $\alpha ,\beta\in {\C}$, thus $<a>$ is isomorphic to ${\C}^2$.   
Let $M_2 ({\C})$ be the unital complex algebra of 
$2\times  2$-matrices with complex entries.  Pick a number $h\notin 
\{0 ,1\}$,  and define the $2\times 2$-matrix  
\[ \pi (a)=\left( \begin{array}{cc}  
0 & h^{-1} \\ h & 0 \end{array}\right) 
\] 
We have $\pi (a)^2 =I$, so by the universal property of $<a>$ we get 
a  representation $\pi :{\C}^2\rightarrow M_2 ({\C})$.  There 
are  no nontrivial (self-adjoint) projections intertwining it.  
However,  $\pi$ is not irreducible (in our stronger sense) as the  
(non-self-adjoint) idempotent $\frac{1}{2}(I-\pi (a))$ is nontrivial  
and intertwines $\pi$.\end{exa} 
\begin{lem} \label{pos} Let $\pi :{\C}(X)\rightarrow B(K)$ be a 
representation, and let $\chi _A\in {\C}(X)$ be the characteristic 
function on $A\subset X$. 
If $e=\pi (\chi_A )$ is different from zero and the identity, then 
the spectrum $\sigma (e)$ of $e$ equals $\{0,1\}$. If $f\geq 0$  
for $f\in {\C}(X)$, 
then $\sigma (\pi (f))\subset [0,\infty )$.  
\end{lem} 
\begin{proof} Suppose that $\lambda\notin \{0,1\}$. Then 
$\chi_A -\lambda I$ is invertible  
in ${\C}(X)$, so $e-\lambda I$ is invertible, 
and thus $\lambda\notin\sigma (e)$, showing that $\sigma 
(e)\subset\{0,1\}$. As $e\neq 0$ pick a nonzero vector $\xi\in 
eK\subset K$. Then $e\xi =\xi$, so $1\in\sigma (e)$. And similar, as 
$e\neq I$, pick a nonzero 
$\xi\in (I-e)K$. Then $e\xi =0$, so $0\in\sigma (e)$, which proves 
the first assertion. 
For the second assertion suppose first 
that $a>0$ and $-a\in\sigma (\pi(f))$. Then, by definition of 
spectrum, $\pi (f+aI)=\pi (f)+aI$ is not 
invertible. But $f+aI\in {\C}(X)$ is strictly positive, and 
therefore invertible, showing that $\pi (f+aI)$ is invertible, a 
contradiction. Therefore there are no (strictly) negative numbers in 
$\sigma (\pi (f))$. 
Next suppose that the imaginary part of $a\in\sigma (\pi (f))$ is 
nonzero.     Then we may write $a+1=\rho\exp (i\theta )$ for numbers 
$\rho ,\theta$ 
with $\rho >0$ 
and $\sin \theta\neq 0$. Thus by spectral calculus \cite{Ru} 
$-4\sin^2\theta <0$ belongs to 
$$\sigma (\pi (((f+I)\rho^{-1}-(f+I)^{-1}\rho )^2 ))$$ 
and clearly 
$$((f+I)\rho^{-1}-(f+I)^{-1}\rho )^2 \geq 0,$$ 
which is impossible according to the previous part of the proof. 
Therefore the imaginary part of $a\in\sigma (\pi (f))$ cannot  be 
nonzero, nor can $a$ be negative, so $a\geq 0$, as desired. 
\end{proof} 
Recall that a filter over a set $X$ is a collection $F$ of 
nonempty subsets of $X$ containing $X$, closed under finite 
intersections, and such that $V\in F$ whenever $U\in F$ and $U\subset 
V\subset X$. 
A filter $F$ is called an ultrafilter if for every $U\in X$ either 
$U\in F$ or its complement $U^c\in F$. This happens if and only if 
$F$ is a maximal filter.  
A filter $F$ is called $\sigma$-complete if it is closed under 
countable  intersections. More generally, if $\kappa$ is any  regular
uncountable cardinal, 
cf. \cite{TJ}, we say that a filter $F$ is $\kappa$-complete if it is 
closed under intersections of less than $\kappa$ sets. 
It is known that an ultrafilter $F$ is $\sigma$-complete if and only
if there is no countable partition $\{U_n \}_{n=0}^{\infty}$ of $X$ 
such that $U_n\notin F$ for all $n$. In fact, a similar result holds 
for $\kappa$-complete ultrafilters, see Exercise 27.2 in \cite{TJ}, 
namely an ultrafilter is $\kappa$-complete if and only if for any 
partition of $X$ consisting of less than $\kappa$ elements, there 
exists a unique element belonging to the ultrafilter.
\bigskip 
We will now see how an irreducible representation of $B$ gives  
a $\sigma$-complete ultrafilter over $X$. 
Let $U\subset X$ and define the element $\chi_U\in B$ by  
$\chi_U (x)=I_x$ 
if 
$x\in U$ and $\chi_U (x)=0$ if $x\in U^c$, where $I_x$ is the 
identity element of $B(K_x )$.  
Let $\pi :B\rightarrow B(K)$ be an irreducible representation of $B$ 
on a Banach space $K$, 
and define $F_{\pi}$ to be the collection of subsets of $X$ given by 
$$F_{\pi}=\{U\subset X\ |\ \pi (\chi_U )=I_{\pi}\},$$ 
where $I_{\pi}$ is the unit in $B(K)$. 
\begin{lem}\label{ultra} Let notation be as above. The collection 
$F_{\pi }$  is a $\sigma$-complete ultrafilter over $X$.
\end{lem} 
\begin{proof} First note that $\pi (0)=0$ and  
$\pi (\chi_X )=\pi (I)=I_{\pi}$ assure that  $F_{\pi}$ is a 
collection of nonempty subsets of $X$ containing $X$. 
If $U,V\in F_{\pi}$, then $U\cap V\in F_{\pi}$, because 
$\chi_{U\cap V}=\chi_U\chi_V$, so 
$\pi (\chi_{U\cap V})=\pi (\chi_U )\pi (\chi_V )=I_{\pi}^2 =I_{\pi}$. 
Let $U\subset X$. 
Note that $\chi_U b=b\chi_U$, so $\pi (\chi_U )\pi (b)=\pi (b)\pi 
(\chi_U )$ holds for all $b\in B$, and furthermore 
$\pi (\chi_U )^2 =\pi (\chi_U )$ as $\chi_U ^2 =\chi_U$. Therefore  
$\pi (\chi_U )$ being an idempotent intertwining the (strongly)  
irreducible $\pi$, 
has to be either $0$ or $I_{\pi}$. 
Let $U\in F_{\pi}$ and suppose that $U\subset V$. Then 
$V\in F_{\pi}$, because 
$\chi_V -\chi_U =\chi_{V\backslash U}$ and so 
$$\pi (\chi _V )-I_{\pi}=\pi (\chi_{V\backslash U})\in\{ 
0,I_{\pi}\},$$ having $\pi (\chi_V )=I_{\pi}$ as the only solution  
within $\{0,I_{\pi}\}$. Thus $V\in F_{\pi}$ 
as well. 
If $I:x\mapsto I_x$, $x\in X$, is the unit of $B$, then 
$I-\chi_U =\chi_{U^c}$, and thus $I_{\pi}-\pi (\chi_U )=\pi 
(\chi_{U^c })$ for any $U\subset X$. Hence $U\in F_{\pi}$ or $U^c\in 
F_{\pi}$, because if $U\notin F_{\pi}$, then $\pi (\chi_U )=0$, 
saying that 
$\pi (\chi_{U^c})=I_{\pi}-\pi (\chi_{U})=I_{\pi}$, and so 
$U^c\in F_{\pi}$. 
We have thus shown that $F_{\pi}$ is an ultrafilter over $X$.  
It remains to 
show that it is in fact $\sigma$-complete.   
Suppose by ad absurdum that we have a countable (disjoint) partition  
$\{U_n\}_{n=0}^{\infty}$ of $X$ such that 
$U_n\notin F_{\pi}$, i.e., $\pi (\chi_{U_n})=0$ for all $n$. 
Define $b^{\infty},b^n\in B$ by 
$b^{\infty}(x)=nI_x$ for $x\in U_n$, $n\in\{0,...,\infty\}$, and $b^n 
(x)=kI_x$ for $x\in U_k$, $k\in\{0,...,n\}$, whereas 
$b^n (x)=nI_x$ for $x\in U_k$ and $k\in\{n+1,...,\infty\}$. 
Then 
$$b^n -nI=\sum_{k=0}^{n-1}(k-n)\chi_{U_k},$$ 
and so 
$$\pi (b^n -nI)=\sum_{k=0}^{n-1}(k-n)\pi (\chi_{U_k}) 
=\sum_{k=0}^{n-1}(k-n)0=0,$$  
which shows that 
$\pi (b^n )=\pi (nI)=nI_{\pi}$ 
for all $n$. 
But $b^{\infty}(x)=\alpha (x)I_x$ and $b^n (x)=\beta (x)I_x$ with  
 $\alpha,\beta\in{\C}(X)$. Hence we may identify $b^{\infty}(x)$  
 with $\alpha$ and $b^n (x)$ with $\beta$,   
and using spectral calculus and Lemma \ref{pos}, we get for all $n$ 
that $$\sigma (\pi (b^{\infty} ))-n=\sigma (\pi (b^{\infty}-b^n )),$$ 
which is impossible, as the spectrum of any element in a Banach 
algebra is nonempty \cite{Ru}, and as Lemma \ref{pos} says, both 
$\sigma (\pi (b^{\infty}))$  and $\sigma (\pi (b^{\infty}-b^n ))$ 
consist solely of positive numbers.  \end{proof} 
\begin{rem}\label{rem:unbdd} The previous Lemma is the key step which 
allows us to pass from ultrafilters to $\sigma$-complete 
ultrafilters, hence to $\gamma$-complete ones, $\gamma$ being the 
least measurable cardinal, as shown below. It is also the only place 
where the existence of unbounded elements in $B$ is used. Things do 
not change if we consider a Hilbert space $H$ and allow unbounded 
elements of $B$ to be represented by unbounded closed linear 
operators on $H$. Indeed, in this case previous Lemma holds the same, 
hence again irreducible representations gives rise to 
$\gamma$-complete ultrafilters. \end{rem}
Recall that a filter $F$ over $X$ is principal if there is a 
$V\subset X$  such that 
$$F=\{U\subset X\ |\ V\subset U\}.$$ 
In case $F$ is an ultrafilter, it is maximal, so $V$ has to be a one 
point set. 
An uncountable set $X$ is of measurable cardinality if there exists a 
nonprincipal $X$-complete ultrafilter over it, cf. p. 297 in 
\cite{TJ}. 
In particular the least cardinal (if any!) on which there is a 
$\sigma$-complete nonprincipal ultrafilter is measurable.
We identify here the set $X$ with its cardinality. We shall say that
a set has less than measurable cardinality if its cardinality is 
smaller than the first measurable cardinal.
\begin{lem}\label{lem2a} Let $F$ be a $\sigma$-complete ultrafilter 
over $X$. If $\gamma$ is the least measurable cardinal, then $F$ is a 
$\gamma$-complete ultrafilter over $X$. In particular, if $X$ has 
less than measurable cardinality, $F$ is principal. \end{lem} 
\begin{proof} Let $\{X_i\}_{i\in\kappa }$ be a partition of
$X$, $\kappa<\gamma$.  We must show that there exists $j\in\kappa$ 
such that $X_j\in F$.  Define a collection of subsets of $\kappa$ by 
$$F_{\kappa}=\{I\subset\kappa\ |\ \cup_{i\in I}X_i\in F\}.$$ 
We claim that $F_{\kappa}$ is a $\sigma$-complete ultrafilter over 
$\kappa$: Clearly $X\in F$ implies that $\kappa\in F_{\kappa}$, and 
$F_{\kappa}$ consists of nonempty subsets only.
If $I,J\in F_{\kappa}$,  
then as $\{X_i\}_{i\in\kappa}$ is a partition of $X$, we have 
$$\cup_{i\in I\cap J}X_i =\cup_{i\in I}X_i\cap\cup_{i\in J}X_i\in 
F,$$ so $I\cap J\in F_{\kappa}$. 
If $J\in F_{\kappa}$ and $I\supset J$, then $I\in F_{\kappa}$ as 
$\cup_{i\in I}X_i\supset\cup_{i\in J}X_i\in F$.  
Furthermore, if $J\notin F_{\kappa}$ then $J^c\in F_{\kappa}$, 
because  $$\cup_{i\in J^c }X_i =(\cup_{i\in J}X_i )^c$$  
as $\{X_i\}_{i\in\kappa}$ is a partition.  
Finally, if $\{\kappa_n\}_{n=0}^{\infty}$ is a countable partition  
of $\kappa$, then $\{Y_n\}_{n=0}^{\infty}$ with  
$Y_n =\cup_{i\in\kappa_n }X_i$ is a countable partition of $X$. Thus 
there exists a number $m$ such that $\cup_{i\in\kappa_m}X_i\in F$, 
i.e., 
$\kappa_m\in F_{\kappa}$, so $F_{\kappa}$ is a $\sigma$-complete 
ultrafilter  over $\kappa$.  
We conclude that $F_{\kappa}$ is a principal filter. Indeed, 
according  
to Lemma 27.1 in \cite{TJ}, the least cardinality for which 
there exists a nonprincipal $\sigma$-complete ultrafilter is 
measurable, therefore $F_{\kappa}$ is principal, hence
there exists $j\in\kappa$ such that $X_j\in F$ as desired. 
\end{proof} 
Now suppose that $0<\dim K_x <\infty$ for all $x\in X$. 
We will see how a $\sigma$-complete ultrafilter $F$ over $X$ together 
with  a collection of bases, one for each $K_x$, $x\in X$, give  
rise to a finite dimensional  
irreducible representation $\pi_{F}$ of $B$.  
Define for any natural number $n$ the set  
$$\Omega_n =\{x\in X\ |\ \dim K_x =n\}\subset X.$$ Clearly, the 
collection $\{\Omega _n \}_{n=0}^{\infty}$  is a countable partition 
of $X$, 
and thus there exists a natural number $n(F)$ such that 
$\Omega_{n(F)}\in F$. This number $n(F)$ is also unique, since if any 
other member of the partition  would belong to $F$, so would their 
intersection, which is empty, a  contradiction. For any $M\in B ({\C}^m )$ 
we denote by $M_{ij}$  
the matrix coefficients  
of $M$ with respect to the standard  
basis on ${\C}^m$. Let $b\in B$ and denote 
by $b(x)_{ij}$ the entries of the matrix $b(x)$  
with respect to the chosen basis of $K_x$.   
Define $I_M\subset X$ by 
$$I_M =\{x\in\Omega_{n(F)}\ |\ b(x)_{ij}=M_{ij}\}.$$ 
Note that $\Omega_{n(F)}^c$ (which does not belong to $F$) together  
with the collection $\{I_M\ |\ M\in B({\C}^{n(F)})\}$  
form a partition of $X$ into as many parts as there 
are real numbers. 
We now use the fact that a $\sigma$-complete ultrafilter is also   
$\gamma$-complete, where $\gamma$ is the least
measurable cardinal, see Lemma \ref{lem2a}. 
So in particular, for the partition above we know that  there is a 
unique  
$M(b)\in B({\C}^{n(F)})$ such that $I_{M(b)}\in F$. Define  
$$\pi_F (b)=M(b),$$   
and so we get a map $\pi_F :B\rightarrow B({\C}^{n(F)})$. We show 
that $\pi_F$ is indeed an irreducible  representation of $B$ 
on ${\C}^{n(F)}$:  
In proving this let's first agree on considering the operators  
as matrices with respect to the given bases, 
to avoid confusion with identifications. 
Then note that if $\Omega\in F$ and $b(x)=b(y)$ (as matrices) for all 
$x,y\in\Omega$, then $\pi_F (b)=b(x)$ for all $x\in\Omega$, because 
$\pi_F (b)=M(b)=b(x)$ for $x\in\Omega\cap I_{M(b)}\neq\emptyset$. 
Consider now $b_1 ,b_2\in B$ and form $I_{M(b_1 )}$ and 
$I_{M(b_2 )}$ in $F$ as prescribed above. As $F$ is a filter,  
$\Omega =I_{M(b_1 )}\cap I_{M(b_2 )}\in F$, so $b_i (x)=b_i (y)$ for 
all  $x,y\in\Omega$. Thus for $x\in\Omega$ it follows from the 
previous argument that 
$$\pi_F (b_1 b_2 )=(b_1 b_2 )(x)=b_1 (x)b_2 (x)=\pi_F (b_1 )\pi (b_2 
).$$ Similarly one proves that $\pi$ is additive, and it is clearly 
unital.  Finally, $\pi (B)=B({\C}^{n(F)})$, indeed if $M\in 
B({\C}^{n(F)})$, define $b\in B$ by 
\[ 
b(x)=\left\{ 
\begin{array}{cc}  
		M & x\in\Omega_{n(F)}  \\ 
		0 & {\mathrm otherwise.} 
\end{array} 
	\right. 
\] 
Clearly $\pi_{F}(b)=M$, so $\pi_F$ is an irreducible 
representation.     In fact, we have the following result.  
\begin{thm} \label{prop1} Let $X$ be an arbitrary set. Suppose that 
the spaces $K_x$ are nonzero and finite dimensional for all $x\in X$. 
The set of equivalence classes of irreducible representations of $B$ 
is in 1-1 correspondence with the set of $\sigma$-complete 
ultrafilters on $X$. More precisely, denote by $[\pi ]$ the 
equivalence class of an
irreducible representation $\pi$ of $B$ on a Banach space $K$, and 
let  $\cal E$ be the set of all such equivalence classes (for all 
Banach spaces $K$). Denote by $\cal F$ the set of all 
$\sigma$-complete ultrafilters over $X$. Then the assignments 
$$[\pi ]\mapsto F_{\pi }$$ 
from $\cal E$ to $\cal F$ and  
$$F\mapsto [\pi_F ]$$ 
from $\cal F$ to $\cal E$, where $F_{\pi}$ is defined as in Lemma 
\ref{ultra} and $\pi_F$ as above, are well defined, and they are 
inverses to each other. \end{thm} 
\begin{proof} The first assignment is well defined, 
as $F_{\pi}=F_{\pi '}$ for two equivalent representations $\pi$ and 
$\pi '$ 
of $B$.  
As for the second assignment notice that $\pi_F$ is defined 
using $F$ and a chosen collection of bases for each $K_x$. To prove 
that 
it is well defined, we therefore need to prove that $[\pi_F ]$ is 
independent 
of this choice of bases. Thus say that we have two collections $(1)$  
and $(2)$ of bases for the spaces $K_x$ and form the two associated 
irreducible representations $\pi_F^1$ and $\pi_F^2$  
of $B$ on ${\C}^{n(F)}$, respectively. We need to prove that they 
are    equivalent.  
As $(1)$ and $(2)$ are collections of bases, we know that to each 
$x\in X$ there exists an invertible matrix $S(x)\in B(K_x )$ taking 
the basis 
in $K_x$ from $(1)$ to the basis in $K_x$ from $(2)$. 
Define for any (finite dimensional complex) square matrix $S$ the 
set  $$X_S =\{ x\in X\ |\ S(x)=S\}.$$  
The collection $\{ X_S\}$ of subsets of $X$, where the index $S$ 
ranges  over all quadratic matrices 
$S\in\cup_{n=0}^{\infty}{\bf M}_n ({\C})$, is clearly a partition  
of $X$ into no more parts than there are real numbers. As $F$ is 
a $\sigma$-complete ultrafilter, it follows from Lemma \ref{lem2a} 
that there exists a square matrix $S(F)$ such that $X_{S(F)}\in F$. 
Let $b\in B$ and consider the  
operators $\pi_F^i (B)\in B({\C}^{n(F)})$, $i=1$ or $i=2$, as    
matrices with respect to the standard bases (1) and (2). Furthermore, 
let $I_{M(b)}^i\in F$ be as in the definition of the representations 
$\pi_F^i$. Then as $F$ is a filter, we have 
$$\Omega =X_{S(F)}\cap I_{M(b)}^1\cap I_{M(b)}^2\in F.$$    
As we noted just before Theorem \ref{prop1}, we thus get for  
any $x\in\Omega$, 
with $b_i (x)$ the matrices with respect to the bases  
for $K_x$ from (i), that  
$$\pi_F^2 (b)=b_2 (x)=S(x)b_1 (x)S(x)^{-1}=S(F)b_1 (x)S(F)^{-1} 
=S(F)\pi_F^1 (b)S(F)^{-1}.$$ 
Hence we have proved that $\pi_F^1$ and $\pi_F^2$ are equivalent, and 
thus  that the second assignment in the theorem is well defined. 
We now prove that these assignments are inverses to each other.  
If we start with  
$F\in {\cal F}$, form $[\pi_F ]\in {\cal E}$ and   
$F_{\pi_F}\in{\cal F}$, then we clearly end up with $F$. 
Thus we only need to see that going in the other direction gives us 
the identity map on ${\cal E}$. 
Let $[\pi ]\in {\cal E}$, and form 
$F_{\pi}\in {\cal F}$ as prescribed;  $$F_{\pi}=\{U\subset X\ |\ \pi 
(\chi _U )=I_{\pi }\},$$ 
where now $\pi$ is an irreducible representation of $B$ on $K$. 
To construct the irreducible representation $\pi_{F_{\pi}}$, we 
choose 
and fix a collection of bases for all the spaces $K_x$. Denote by 
$n(\pi )$ the natural number $n(F_{\pi})$ associated to the 
$\sigma$-complete ultrafilter $F_{\pi}$ such that $\Omega_{n(\pi 
)}\in F_{\pi}$. 
Define $B_{n(\pi )}$ by  
$$B_{n(\pi )}=\prod_{x\in\Omega_{n(\pi )}}B(K_x ).$$ 
We complete the proof by showing that 
$[\pi_{F_{\pi}}]=[\pi]$, and we proceed in five steps: 
First 1) we claim that   
$$\pi (B)=\pi (B_{n(\pi )}).$$ 
To see this notice that for $b\in B$, we have $$b=\chi_{\Omega_{n(\pi 
)}}b+(1-\chi_{\Omega_{n(\pi )}})b,$$ thus  
$$\pi (b)=\pi (\chi_{\Omega_{n(\pi )}}b)$$ 
and  
$$\chi_{\Omega_{n(\pi )}}b\in B_{n(\pi )}.$$ 
Secondly 2), obviously 
$B({\C}^{n(\pi )})$ equals  
$$C=\{b\in B_{n(\pi )}\ |\ \exists c,\forall x,\ b(x)_{ij} 
=c_{ij}\},$$   when the operators are considered as matrices with 
respect to the bases for  the spaces $K_x$ and ${\C}^{n(\pi )}$. 
Thirdly 3) we have  
$$\pi (B_{n(\pi )})=\pi (C).$$ 
For let $b\in B_{n(\pi )}$, define for $m\in B({\C}^{n(\pi )})$ 
with respect to bases for $K_x$ and ${\C}^{n(\pi )}$, the set    
$$X_m =\{x\in\Omega_{n(\pi )}\ |\ b(x)_{ij} =m_{ij}\},$$ 
and consider the partition  
$$\{X_m\ |\ m\in B({\C}^{n(\pi )})\}$$ 
of $\Omega_{n(\pi )}$. By adding $\Omega_{n(\pi )}^c\notin F_{\pi}$ 
to this  partition, we get a partition of $X$ into no more parts than 
there are real   numbers. It follows then from Lemma \ref{lem2a}, 
that there exists  
a unique $m(b)$ such that $X_{m(b)}\in F_{\pi}$. Then $\pi (b)=\pi 
(m(b))$, because 
$$b-m(b)=(b-m(b))\chi_{X_{m(b)}}+(b-m(b))(1-\chi_{X_{m(b)}})$$ 
$$=(b-m(b))\chi_{X_{m(b)}^c}.$$ 
Now $X_{m(b)}^c\notin F_{\pi }$, so  
$$\pi (b-m(b))=\pi (b-m(b))\pi (\chi_{X_{m(b)}^c})=0.$$ 
>From 1),2) and 3) above we have $\pi (B)=\pi (B({\C}^{n(\pi )}))$. 
By nontriviality of $\pi$, we get that $\pi (B({\C}^{n(\pi )}))$ 
is isomorphic to $B({\C}^{n(\pi )})$, as otherwise $\ker \pi$ 
would be a nontrivial ideal in the simple algebra $B({\C}^{n(\pi 
)}))$. Thus $\pi (B)$ is isomorphic to the algebra $B({\C}^{n(\pi 
)}))$,  and as $\pi$ is irreducible, $K$ is isomorphic to ${\C}^{n(\pi )}$. 
We can now add the final step of the proof. By definition of   
the representation $\pi_{F_{\pi}}$ we have for $b\in B$, that 
$$\pi_{F_{\pi}}(b)=M(b),$$ 
where $$\{x\in\Omega_{n(\pi )}\ |\ b(x)_{ij} =M(b)_{ij}\}\in 
F_{\pi}.$$ 
But from above, we see that  
$$\pi (b)=\pi (M(b))=M(b),$$ 
by the identification of $\pi (B)$ with $B({\C}^{n(\pi )})$, 
so $\pi (b)=\pi_{F_{\pi}}(b)$ and $[\pi ]=[\pi_{F_{\pi}}]$ as 
wanted.  \end{proof} 
\begin{rem} We may obviously generalize this results to the case 
where the spaces $K_x$ are Hilbert spaces with less than measurable 
dimension. More precisely, there is a 1-1 correspondence between 
equivalence classes of irreducible representations of $B$ (= 
algebraic product of the $B(K_x )$'s for $x\in X$) and 
$\sigma$-complete ultrafilters over
 $X$.
\end{rem}
\section{Characterization of Representations on $B$}  
Our main result in this section is Theorem \ref{Th2}. It 
says that whenever the index set the product in $B$ is taken over, 
is of less than measurable cardinality, then the only 
representations of $B$ on Banach spaces, are finite 
sums of those one gets by projecting down on the factors in the 
product. 
\begin{thm}\label{Th1} Let $X$ be a set of less than measurable 
cardinality, and let $\pi$ be a representation of $B$ on a Banach 
space $K$. Then $\pi$ is determined by the representations of the 
algebras $B(K_x )$ for a finite number of $x$'s. More precisely, 
define for $Z\subset X$ the subalgebra $B_Z =\prod_{x\in Z}B(K_x)$ of 
$B$, and let $\pi_
Z$ be the representation of $B$ on $K$ obtained by projecting down 
$B\rightarrow B_Z\subset B$ and then restricting $\pi$ to $B_Z$: 
$\pi_Z(b)=\pi(\chi_{Z}b)$. Define $Y\subset X$ by 
$$Y=\{x\in X\ |\ \pi_x\neq 0\}.$$  
Then 1) Y is a finite set, and 2) $\pi_Y$ is equivalent to $\pi$.  
\end{thm} 
\begin{proof} Define a collection $U_{\pi}$ of subsets of $X$ 
by $$U_{\pi}=\{U\subset X\ |\ \pi (\chi_U )\neq 0\}.$$  
The collection $U_{\pi }$ has the inclusion property, because 
say that we have a set $A\in U_{\pi}$ and a set $B$ with  
$A\cap B=\emptyset$ such that $A\cup B\notin U_{\pi}$. Then  $0=\pi 
(\chi_{A\cup B})=\pi (\chi_A )+\pi (\chi_B )$, so $\pi (\chi_B )=-\pi 
(\chi_A )$, which is impossible as both $\pi (\chi_A )\neq 0$ and 
$\pi (\chi_B )$ are idempotents. 
But $U_{\pi}$ is not necessarily an ultrafilter (which holds only 
when $\pi$ is irreducible).  
However, there exists a finite maximal  
partition $\{X_i\}_{i=1}^n$ of $X$ with the property  
that $\pi (\chi_{X_i })\neq 0$ for all $i$, and furthermore,  
such that the collections 
$$F_i =\{U\in U_{\pi }\ |\ U\subset X_i\}$$ 
are all $\sigma$-complete ultrafilters over $X_i$.  
Assume for the moment that we have such a partition of $X$.  
Observe first that $U\in U_{\pi}$ if and only if there exists $i$ such 
that $U\cap X_i\in F_i$. 
As $X$ is a set of less than measurable cardinality, 
so are the sets $X_i$, and so by 
Lemma \ref{lem2a} all filters $F_i$ are based on (unique) points 
$x_i\in X_i$. Thus $U\in U_{\pi }$ if and only if there exists $i$ 
such that $x_i\in U$. 
Now set $W=\{x_i\}_{i=1}^n$.  Certainly $W$ is a finite subset of 
$X$. 
As $W^c\cap X_i =X_i\backslash \{x_i \}\notin U_{\pi}$ and $X_i\in 
U_{\pi}$, it follows by the intersection property of $U_{\pi}$ that  
$W^c\notin U_{\pi}$, so $\pi (\chi_{W^c})=0$. For $b\in B$ we thus 
get $$\pi (b)=\pi ((\chi_{W}+\chi_{W^c})b)=\pi (\chi_{W}b)=\pi_W 
(b),$$ so $\pi$ is equal (and so a fortiori equivalent) to $\pi_W$ 
and $W=Y$. 
Therefore we are done if we can prove that a  
finite partition $\{X_i\}_{i=1}^n$  
of $X$ as described above does indeed exist. 
We first prove that if  
\begin{equation} 
\label{partition1} 
\exists\ a\ finite\ maximal\ partition\ \{X_i\}_{i=1}^n \ of\ X\ 
such\ that\ \pi (\chi_{X_i})\neq 0, 
\end{equation} 
then the collections $F_i$ are all $\sigma$-complete ultrafilters  
over $X_i$. 
To this end we first note that if $U\subset X_j$  
with $\pi (\chi_U )\neq 0$, then $\pi (\chi_U )=\pi (\chi_{X_j })$. 
This follows because 
$$\pi (\chi_{X_j} -\chi_U )=\pi (\chi_{X_j\cap U^c })=0,$$ 
as otherwise $\{X_i\}_{i\neq j}$, $U$ and $X_j \cap U$ would be a 
partition of $X$ satisfying (\ref{partition1})  
and which is strictly larger than the maximal one.  
Hence for any $U\subset X_i$, we either have $\pi (\chi_U )=0$ 
or $\pi (\chi_U )=\pi (\chi_{X_i} )$. Now proceed as in the proof of 
Lemma \ref{ultra} with $I_{\pi}$ replaced by $\pi (\chi_{X_i})$, to 
prove that the collections $F_i$ are all $\sigma$-complete 
ultrafilters over $X_i$. 
Hence we are left with proving the existence of a finite partition  
$\{X_i\}_{i=1}^n$ with property (\ref{partition1}). As $\pi (\chi_X 
)\neq 0$,  this is clearly equivalent to saying that 
$$\sup \{n\in {\N}\ |\ \exists\ a\ partition\ \{X_i \}_{i=1}^n\ 
of\ X\ such\ that\ \pi (\chi_{X_i })\neq 0\}<\infty .$$ 
Using the inclusion property of $U_{\pi}$, one sees that the  
negation of this statement is equivalent to saying that  
\begin{equation} 
\label{partition2} 
\forall\ n\in {\N}\ \exists\ a\ partition\ \{X_i\}_{i=1}^n\ of\ X\ 
with\  \pi (\chi_{X_i })\neq 0,\ \forall\ i.  
\end{equation} 
Going for ad absurdum we may complete the proof by showing that 
property (\ref{partition2}) leads to a contradiction.  We will obtain 
this contradiction by constructing a countable partition $\{ 
X_i\}_{i=1}^{\infty}$ of $X$ such that $\pi (\chi_{X_i })\neq 0$.  
This is sufficient, because say we have such a partition 
$\{X_i\}_{i=1}^{\infty}$ of $X$.  By Lemma \ref{pos} we conclude that 
$1\in \sigma (\pi (\chi_{X_i}))$ for all $i$.  Define $f\in B$ by 
$f=\sum_{n=1}^{\infty}n\chi_{X_n}$ and note that $$\sigma (\pi 
(f))=\bigcup_{n=1}^{\infty}n(\sigma (\pi (\chi_{X_n}))).$$ Thus 
$\sigma (\pi (f))={\N}$, which is a contradiction as the spectrum 
of an operator in a Banach space $B(K)$ is bounded.  Thus assume that 
property (\ref{partition2}) holds for $X$.  Now to construct the 
desired countable partition of $X$, we rely on the following 
property: 
Namely, if $U_1$ and $U_2$ form a partition of a set $U$ having 
property (\ref{partition2}), then $U_1$ or $U_2$ have property 
(\ref{partition2}).  Assume by absurdum that neither $U_1$ nor $U_2$ 
have property (\ref{partition2}).  Then there exist numbers $n_1,n_2$ 
such that $\{U_j^i\}_{i=1}^{n_j}$ are maximal partitions of $U_j$ 
such 
that $\pi (\chi_{U_j^i})\neq 0$ for all $i=1,\dots,n,j=1,2$.  Let $N$ 
be $\max(n_1,n_2)$.  Thus for $n>N$ there are no partitions 
$\{U_j^i\}_{i=1}^n$ of $U_j$ such that $\pi (\chi_{U_j^i})\neq 0$ for 
all $i=1\dots n,j=1,2$.  Since $U$ has property (\ref{partition2}), 
we 
may choose a partition $\{X_i\}_{i=1}^{2N+1}$ of $U$ with $\pi 
(\chi_{X_i })\neq 0$.  Then define the sets $U_j^i\subset U_j$ by 
$U_j^i =U_j\cap X_i$.  They give partitions of $U_1,U_2$, 
respectively, and, for any $i=1,\dots,2N+1$ there exists $j$, 
depending on $i$, such that $\pi(\chi_{U_j^i})\ne0$, namely there 
exists $j=1,2$ such that at least $N+1$ elements of the partition 
$\bigcup_{i=1}^{2N+1}U_j^i$ have the property that 
$\pi(\chi_{U_j^i})\ne0$.  By adding the complement to one of the 
elements of the partition and using the inclusion property for 
$U_{\pi}$, we therefore obtain a partition $\{\Omega_i\}_{i=1}^{N+1}$ 
of one of the sets $U_j$ such that $\pi (\chi_{\Omega_i})\neq 0$ for 
all $i$, a contradiction with the above.  Thus either $U_1$ or $U_2$ 
have property (\ref{partition2}).  Consider now again the set $X$.  
It 
has property (\ref{partition2}), so split it up in two parts, say 
$X_1^0$ and $X_2^0$ with $\pi (\chi_{X_i^0})\neq 0$.  Then by the 
above, (at least) one of them have property (\ref{partition2}), say 
$X_1^0$, so we may split this one up in two parts $X_1^1$ and $X_2^1$ 
with $\pi (\chi_{X_i^1})\neq 0$.  Again one of these must have 
property (\ref{partition2}), say $X_1^1$, so we may split it up in 
two 
parts $X_1^2$ and $X_2^2$ with $\pi (\chi_{X_i^2})\neq 0$, and one of 
them must have property (\ref{partition2}) again, say $X_1^2$.  This 
way we construct an infinite sequence of disjoint subsets $\{ 
X_1^i\}_{i=1}^{\infty}$ of $X$ such that $\pi (\chi_{X_2^i})\neq 0$.
Now set  $$X_0=(\bigcup_{i=1}^{\infty}(X_2^i ))^c.$$ 
Clearly $X_{0}\supset X_{2}^0$, so, by the inclusion property of 
$U_{\pi}$,  we have that $\pi (\chi_{X_1 })\neq 0$. 
Define for $i\geq 1$ the sets $X_i =X_2^i$. 
Then clearly the collection $\{X_i\}_{i=0}^{\infty }$ 
is a countable partition of $X$ such that $\pi (\chi_{X_i})\neq 0$ 
for all $i$. This concludes the proof.  
\end{proof} 
We give the following lemma without a proof, and refer to standard   
texts on $C^*$-algebras, cf. \cite{FD}. 
\begin{lem}\label{lemmaA} Let $A$ be a finite dimensional 
$C^*$-algebra, so $A$ is isomorphic to a finite sum of full matrix 
algebras ${\bf M}_{k(i)}({\C})$; $$A\simeq \oplus_{i=1}^m {\bf 
M}_{k(i)}({\C}).$$
Define the representations $\pi_i$ of $A$ on ${\C}^{k(i)}$ by 
$$\pi (\oplus_{j=1}^m m_{k(j)})=m_{k(i)}.$$
If $\pi$ is a representation of $A$ on a Banach space $K$, then there 
exist Banach spaces $K(i)$ such that $\pi$ is (weakly) equivalent to 
the representation $\oplus_{i=1}^m \pi_i\otimes I_{K(i)}$, where 
$I_{K(i)}$ is the identity operator on the Banach space $K(i)$. 
\end{lem} 
Let $X$ be a set of less than measurable cardinality and suppose 
that  
all the Banach spaces $K_x$ in the product algebra $B$ are finite  
dimensional.  Denote by $p_x$ the representation of $B$ on $K_x$  
obtained by projecting down on the factor $B(K_x )$, namely  
$p_{x}(b)=b(y)$.  Then, given a finite subset $Y$ of $X$ and a 
Banach  space valued map $\kappa:x\in Y\mapsto K(x)$ , we may 
consider the  following representation of $B$: 
$$ 
p_{\kappa}=\bigoplus_{x\in Y}p_x\otimes I_{K(x)}, 
$$ 
$I_{K(x)}$ being the identity operator on the space $K(x)$. 

\begin{thm}\label{Th2} Let $X$ be a set of less than measurable 
cardinality and suppose that all the Banach spaces $K_x$ in the 
product algebra $B$ are finite dimensional. Then any representation 
$\pi$ of $B$ on a Banach space $K$ is equivalent to $p_{\kappa}$ for 
some finite subset $Y$ in $X$ and some map $\kappa$ as above. 
\end{thm} 
\begin{proof} The theorem is an immediate consequence of 
Theorem 
\ref{Th1} and Lemma \ref{lemmaA}. Note that the $K(x)$'s in the 
statement may be assumed to be nonzero. 
\end{proof} 
The corollary below is an immediate consequence of Theorem 
\ref{Th2}.  But as the proof of this theorem is rather long, we add a 
simplified proof of the corollary, from which the subsequent 
corollaries are proved as well. 
\begin{corol}\label{projned} Let $X$ be a set of less than measurable 
cardinality and suppose that all the Banach spaces $K_x$ in the 
product algebra $B$ are finite dimensional. Then any irreducible 
representation $\pi$ of $B$ is equivalent to $p_{y}$ for some $y\in 
X$. 
\end{corol} 
\begin{proof} Let $F_{\pi}$ be the $\sigma$-complete 
ultrafilter 
over $X$ given by the representation $\pi$, see Lemma \ref{ultra}.  
According to Lemma \ref{lem2a}, since $X$ has less than measurable 
cardinality, $F_{\pi}$ is principal.
Thus for $X$ having less than measurable cardinality there exists  
$y\in X$ such that 
$$F_{\pi}=\{U\subset X\ |\ y\in U\}.$$ For this $y\in X$, we have by 
definition 
of $F_{\pi}$ that $\pi (I_y )=I_{\pi}$, where we consider the unit 
$I_y$ of $B(K_y )$ as an element of $B$.  
As $I=I_y +I_{y^c }$, we have $\pi (I_{y^c})=0$, and thus for all 
$b\in B$ $$\pi (b)=\pi (bI)=\pi (b(I_y +I_{y^c})) 
=\pi (b(y))+\pi (b)\pi (I_{y^c})=\pi (b(y)),$$ 
where $b(y)\in B(K_y )$ is considered an element of $B$ and 
$I_{y^c}$  denotes the element $\chi_{\{y\}^c}\in B$. 
Thus $b(y)\mapsto\pi (b(y))$ is a finite dimensional  
irreducible representation 
of $B(K_y )$ on $K$.  
But the finite dimensional algebra $B(K_y )$ is isomorphic to a  
matrix algebra, 
and we know, cf. for instance \cite{FD}, that every finite 
dimensional  irreducible representation 
of a matrix algebra is equivalent to the identity representation. 
Therefore there exists an invertible linear operator $S:K_y\rightarrow 
K$ such that $\pi (b(y))=S\circ b(y)\circ S^{-1}$ for all $b\in B$. 
Hence $\pi (b)=S\circ p_y (b)\circ S^{-1}$ for all $b\in B$, saying 
that $\pi$ is equivalent to $p_y$, as desired. \end{proof} 
\begin{rem} Let $X$ be a set of less than measurable cardinality. 
Consider the case where all the spaces $K_x$ are finite dimensional, 
and endow them with a Hilbert space structure. Then since the 
representations $p_x$, $x\in X$, are all $*$-representations of $B$ 
on the Hilbert spaces $K_x$, every irreducible representation of 
$B=\prod_{x\in X} B
(K_x )$ is equivalent to a (finite dimensional irreducible) 
$*$-representation of $B$ on some Hilbert space.
\end{rem}
If we choose the Hilbert spaces in $B=\prod_{x\in X}B(K_x )$ all 
to be 1-di\-men\-sion\-al, then  an  easy consequence of Theorem \ref{prop1} 
and Corollary \ref{projned} 
is the result, due to \cite{BBZ} and mentioned in the introduction, 
that when $X$ is a set of less than measurable 
cardinality all characters on ${\C}(X)$ are Dirac measures. 
\bigskip 
\begin{rem} We conclude this paper by observing that even though the
theory of large cardinals is an active field of research concerning 
important set-theoretical questions, measurable cardinals are beyond 
any possibility of construction. In fact, the following results hold: 
Firstly, the existence of sets with measurable cardinality is 
inconsistent with the axiom of constructibility, cf.  Section 31 in 
\cite{TJ}.  Secondly, measurable cardinals are all inaccessible 
cardinals, cf. Lemma 27.2 in \cite{TJ}, which means for instance, 
that their existence is not provable in ZFC (Zermelo-Fraenkels set 
theoretical axioms and the axiom of choice), cf. Theorem 27 in 
\cite{TJ}. Moreover, using G\"odel's second incompleteness theorem in 
its proof, the same theorem says that it cannot be shown that the 
existence of even inaccessible cardinals is consistent with ZFC. 
Therefore Corollary \ref{projned} holds for any 
``concrete'' set $X$. 
\end{rem}
\begin{ack} We thank Paolo Lipparini for conversations. We also thank the
referee for pointing out our attention to reference \cite{BBZ} \end{ack}
 
\end{document}